# Symmetric Spaces and star-representations III. The Poincaré Disc

P.Bieliavsky*, M.Pevzner *

`pbiel@ulb.ac.be, mpevzner@ulb.ac.be`

Université Libre de Bruxelles, Belgium

**Abstract**

This article is a contribution to the domain of (convergent) deformation quantization of symmetric spaces by use of Lie groups representation theory. We realize the regular representation of $SL(2,\mathbb{R})$ on the space of smooth functions on the Poincaré disc as a sub-representation of $SL(2,\mathbb{R})$ in the Weyl-Moyal star product algebra on $\mathbb{R}^2$. We indicate how it is possible to extend our construction to the general case of a Hermitian symmetric space of tube type.

## 1 Introduction

For the convenience of the reader, we start by recalling basic facts in deformation quantization theory which will be useful in explaining what is done in the present article.

### 1.1 Star products

We first recall the notion of deformation quantization (star product) as introduced in the seminal paper [3]. Star products were originally introduced as an autonomous formulation of quantum mechanics in the framework of classical mechanics, that is without reference to a Hilbert space representation. Roughly speaking, if one thinks of a mechanical system as the data of a symplectic manifold $M$ (the phase space) together with some fixed function $H$ on $M$ (the hamiltonian), one usually understands a quantization of the classical system $(M, H)$ as a correspondence (the quantization map) between functions on $M$ (the classical observables) and operators (the quantum observables) on some Hilbert space (the quantum states of the system). The physical information is encoded in the spectral theory of the algebra of quantum operators. In partic­ular, by reading the composition product of quantum operators at the level of

*Research supported by the Communauté française de Belgique, through an Action de Recherche Concertée de la Direction de la Recherche Scientifique

classical functions on $M$ via the quantization map, one should get an associative product on the classical functions encoding the quantum mechanics. Typically, a star product on a symplectic manifold $M$ can be thought as such a transported composition product on functions on $M$. It turns out that star products can be constructed independently of any Hilbert space data. Now precisely, one has

**Definition 1.1** *Let $(M,\omega)$ be a symplectic manifold and $E = C^\infty(M)$ be the commutative algebra of smooth (complex valued) functions on $M$. A formal deformation quantization (or star product) of $(M,\omega)$ is an associative product law on the space $E[[\nu]]$ of formal power series in the parameter $\nu$ with coefficients in $E$. Denoting this product by*

$$(a,b) \mapsto a\star_\nu b, \qquad a,b \in E[[\nu]]. \tag{1}$$

*one requires furthermore the following properties.*

(i) *The map $\star_\nu$ is $\mathbb{C}[[\nu]]$-bilinear ($\mathbb{C}[[\nu]]$ stands for the formal power series with constant complex coefficients).*

(ii) *If for $u,v \in C^\infty(M) \subset E[[\nu]]$, one denotes*

$$u\star_\nu v =: \sum_{k=0}^\infty \nu^k c_k(u,v), \tag{2}$$

*then the maps $c_k : E \to E$ are bidifferential operators on $E$ such that $c_0(u,v) = uv$ and $c_1(u,v) = \{u,v\}$ where $\{\,,\,\}$ stands for the Poisson bracket on $E$ associated with the symplectic structure $\omega$.*

In other words, a formal deformation quantization of a symplectic (or more generally Poisson) manifold is a formal associative deformation of the usual commutative product on functions in the direction of the Poisson structure.

**Definition 1.2** *Two star products $\star_\nu{}^1$ and $\star_\nu{}^2$ on $M$ are called equivalent if there exists a formal series of differential operators $T = Id + \sum_{k\geq 1} \nu^k T_k$ such that for all $u,v \in C^\infty(M)$, one has $T(u\star_\nu{}^1 v) = Tu\star_{\nu^2} Tv$.*

On every Poisson manifold there exists a star product [16], but the most simple example of a star product is the so-called Weyl-Moyal star product on the symplectic vector space $(M = \mathbb{R}^{2n}, \omega = \sum_{i<j} \Lambda_{ij} dx^i \wedge dx^j)$ (the coefficients $\Lambda_{ij}$'s being constant). The Weyl-Moyal star product on $\mathbb{R}^{2n}$ is defined as the $\mathbb{C}[[\nu]]$-bilinear extension to $C^\infty(\mathbb{R}^{2n})[[\nu]]$ of

$$u \star_\nu^0 v := uv + \sum_{k=1}^\infty \frac{\nu^k}{k!} \Lambda^{i_1 j_1} \ldots \Lambda^{i_k j_k} \frac{\partial^k}{\partial x^{i_1} \ldots \partial x^{i_k}} u \, \frac{\partial^k}{\partial x^{j_1} \ldots \partial x^{j_k}} v. \tag{3}$$

It turns out to be associative. The associative algebra $(C^\infty(\mathbb{R}^{2n})[[\nu]], \star_\nu)$ is called the Weyl-Moyal algebra on $\mathbb{R}^{2n}$. It can be seen as a completion of the universal enveloping algebra of the Heisenberg Lie algebra constructed on $(\mathbb{R}^{2n}, \omega)$. It can also be characterized–and that is what matters here– by its automorphism group (see below).

## 1.2 Symmetries

In presence of symmetry i.e. when one is given a symplectic action of some Lie group $G$ on $(M, \omega)$, one may want our quantization map to be equivariant with respect to this action and some unitary representation of $G$ on the state space. At the level of formal deformation quantization, this amounts to require the $G$-invariance of the star product.

**Definition 1.3** *Let $(M, \omega)$ be a symplectic manifold which a Lie group $G$ acts on in a differentiable way by symplectomorphisms (i.e. denoting the action by $G \times M \to M : (g, x) \mapsto \tau_g(x)$, one has $\tau_g^\star \omega = \omega \quad \forall g \in G$). A star product $\star_\nu$ on $(M, \omega)$ is called $G$-invariant if for all functions $u, v \in C^\infty(M)$ and $g \in G$, one has*

$$\tau_g^\star(u \star_\nu v) = \tau_g^\star u \star_\nu \tau_g^\star v. \tag{4}$$

When the action of $G$ preserves a symplectic affine connection on $M$, the existence of a $G$-invariant star product is guaranteed and these star products are classified [4, 13]. In particular, up to a change of parameter $\nu \mapsto \nu'$, the Weyl-Moyal star product is the only star product on the symplectic vector space $\mathbb{R}^{2n}$ which is invariant under $\mathrm{Sp}(n, \mathbb{R}) \times \mathbb{R}^{2n}$ acting on $\mathbb{R}^{2n}$ by affine transformations with respect to the canonical flat connection. In particular, this very basic example already shows an interplay between deformation quantization and geometry in the sense that a star product generally carries a strong geometrical information.

## 1.3 Convergence

Recovering from the data of a formal star product the spectral theory containing the physical information is generally a very hard problem, especially in presence of curvature. This problem appears to be closely related to the one of finding "convergent" (or "non-formal") star products [13, 14, 20]. Many notions of convergence for star products can be found in the literature, but it's probably fair to say that a convergent star product is a star product which admits an asymptotic representation in some operator algebra. For instance the $C^\star$-algebraic situation has been developed by Rieffel in [20] and in earlier works. However, the geometrical affine information encoded in the deformed product in these cases is always associated to a flat affine connection as in the Weyl-Moyal case. It is therefore challenging to try to attach 'à la Rieffel' some topological operator algebra framework, such as a $C^\star$ or Fréchet, to symmetry invariant star products on (symplectic) symmetric spaces. For solvable symplectic symmetric spaces (i.e. whose transvection group is a solvable Lie group), this has been done in [6]. Later, other solvable Lie group actions have been treated in the same spirit in [7]. These works have led to oscillating integral formulae for symmetry invariant associative products on function spaces defined on solvable symmetric spaces. The associated oscillating kernels have been explicitly determined and it turned out that they entirely encode the affine symmetric symplectic geometry of the space at hand (topologically $\mathbb{R}^{2n}$). For a beautiful symplectic geometrical interpretation of the phase part occurring in such an oscillating kernel see [25]. In a different context than star product theory, non-formal quantization of semisimple symmetric spaces has extensively been studied in the important work of Unterberger (see [23] for a survey) by using techniques of harmonic analysis (see also [24] for an operator algebraic approach).

## 1.4 The question attended here and the results obtained

In the present article, we are concerned with investigating the problem of constructing invariant associative products on certain function algebras on semisimple symplectic symmetric spaces by applying to the semisimple situation the same method as used for the solvable case in [6]. This method essentially relies on the fact that formal star products on a symplectic space $M$ are classified (up to equivalences) by the space of formal powers series with coefficients in the second de Rham cohomology space of $M$ [4]. In particular, when the topology of the space is trivial all the star products are equivalent to one another. If one disposes of a Darboux chart on a $G$-space, one should then be able to relate by an equivalence any $G$-invariant star product to the simplest (non-invariant!) one: the Weyl-Moyal star product. For certain symmetric spaces, the above mentioned equivalence $T$ can be obtained quite explicitly from a combination between fairly old techniques in star representation theory (Lie group representation theory in the framework of deformation quantization [14, 1, 9] and a detailed geometrical study. Non-commutative harmonic analysis really enters in the game when passing from this formal situation to a convergent one.

We now precise the construction. Let $G/K$ be an Hermitian symmetric space of tube type. The first step consists in defining an action by star product derivations of the Lie algebra $\mathfrak{g}$ of $G$ on the Weyl-Moyal algebra on $\mathbb{R}^m$. This representation of $\mathfrak{g}$ has already been considered in [5] and further studied in [8]. It is a $\nu$-deformation of the regular infinitesimal action of $\mathfrak{g}$ on the functions on (a dense open Darboux chart domain of) the causal symmetric space of Cayley type $G/H$ associated to $G/K$. For details, see Section 3.

The second step, performed in the present paper, consists in finding a subspace of the Weyl-Moyal algebra to which the representation restricts as (equivalent to) the regular infinitesimal action of $\mathfrak{g}$ on functions on the tube $G/K$. This can be achieved only when the parameter $\nu$ is assigned some particular rational values. In this case, an injective intertwiner

$$C_c^\infty(G/K) \to \text{ Weyl-Moyal algebra}$$

is constructed via the Plancherel formula for $L^2(G/K, \delta)$ where $\delta$ is some finite dimensional unitary representation of $K$ depending on $\nu$. As a consequence the restriction of the Weyl-Moyal product to the image of the intertwiner defines (at least formally) a $G$-equivariant mapping

$$C_c^\infty(G/K) \times C_c^\infty(G/K) \to \text{ Weyl-Moyal algebra} : (u,v) \mapsto u\star_\nu v.$$

This construction yields an interpolation between the regular representations of $G$ on the function on the Cayley symmetric space $G/H$ on the first hand and on the corresponding tube $G/K$ on the second hand. In particular, one is tempted to make the analogy between the latter and the interpolation found by Neretin in the compact-non-compact duality context [18].

We do not know wether it is possible to choose the intertwiner in a way that the algebra generated by the subspace $C_c^\infty(G/K)$ in the Moyal algebra would still be a function algebra on $G/K$ i.e. be part of the regular representation of $G$ on $C^\infty(G/K)$ (in the solvable situation, it is actually possible [6]). This last question certainly deserves quite a bit of interest, and we hope to be able to give some elements of answer to it in a future work.

Throughout this paper, and in order to avoid technicalities we concentrate ourselves to the "test" case $SL(2,\mathbb{R})/SO(2)$. However, we shall give concrete indications on how to extend our construction to the general case of an arbitrary Hermitian symmetric spaces of tube type $G/K$.

**Acknowledgments.** The proof of Theorem 5.1 we adopted is due to E.P. van den Ban (private communication), we express him here our deep gratitude. We would also like to thank G.van Dijk, J.Faraut, and J.A.Wolf for fruitful discussions.

## 2 Covariant $\star$-products

We start by recalling some basics in Hamiltonian systems and then pass to star products. In this section, $(M, \omega)$ is a symplectic manifold and $\mathfrak{g}$ is a finite dimensional real Lie algebra. One assumes one has a representation of $\mathfrak{g}$ as an algebra of symplectic vector fields on $(M, \omega)$. That is, one has a Lie algebra homomorphism $\mathfrak{g} \to \mathcal{X}(M) : X \mapsto X^*$ ($\mathcal{X}(M)$ stands for the space of smooth vector fields on $M$) such that for all $X$ in $\mathfrak{g}$ one has

$$\mathcal{L}_{X^*}\omega = 0,$$

where $\mathcal{L}$ denotes the Lie derivative. One supposes furthermore that this representation of $\mathfrak{g}$ is strongly Hamiltonian which means that there exists a $\mathbb{R}$-linear map

$$\mathfrak{g} \stackrel{\lambda}{\mapsto} C^\infty(M) : X \mapsto \lambda_X,$$

such that

$$\begin{aligned}(i) & \quad d\lambda_X = i_{X^*}\omega, \\ (ii) & \quad \lambda_{[X,Y]} = \{\lambda_X, \lambda_Y\},\end{aligned}$$

where $\{,\}$ denotes the Poisson structure on $C^\infty(M)$ associated to the symplectic form $\omega$.

**Definition 2.1** *A quadruple $(M, \omega, \mathfrak{g}, \lambda)$ with $(M, \omega), \mathfrak{g}$ and $\lambda$ as above is called a strongly hamiltonian system. The map $\lambda : \mathfrak{g} \mapsto C^\infty(M)$ is called the moment mapping.*

**Example 1.** Coadjoint orbits. Let $M = \mathcal{O} \subset \mathfrak{g}^*$ be a coadjoint orbit of a Lie group $G$ with Lie algebra $\mathfrak{g}$. In this case, we denote by $\mathfrak{g} \to \mathcal{X}(M) : X \mapsto X^*$ the rule which associates to an element $X$ in $\mathfrak{g}$ its *fundamental vector field* on $\mathcal{O}$:

$$X^*_x := \frac{d}{dt}|_0 \operatorname{Ad}^*(\exp(-tX))x,$$

where $\operatorname{Ad}^*(g)x$ denotes the coadjoint action of the element $g \in G$ on $x \in \mathfrak{g}^*$. The formula $\omega^\mathcal{O}_x(X^*, Y^*) := \langle x, [X, Y]\rangle$ (with $X, Y \in \mathfrak{g}$) then defines a symplectic structure called after, Kirillov, Kostant and Souriau, the *KKS symplectic form* on $\mathcal{O}$. Defining, for all $X \in \mathfrak{g}$, the function $\lambda_X \in C^\infty(\mathcal{O})$ by

$$\lambda_X(x) := \langle x, X\rangle,$$

one then has that the quadruple $(\mathcal{O}, \omega^{\mathcal{O}}, \mathfrak{g}, \lambda)$ is a strongly Hamiltonian system. In the setting of deformation quantization there is a natural way to define the equivariant quantization of a classical hamiltonian system [1],[3],[14].

**Definition 2.2** *[1] Let $(M, \omega, \mathfrak{g}, \lambda)$ be a classical strongly hamiltonian system. A star product $\star_\nu$ on $(M, \omega)$ is called $\mathfrak{g}$-covariant if for all $X, Y \in \mathfrak{g}$ one has*

$$\lambda_X \star_\nu \lambda_Y - \lambda_Y \star_\nu \lambda_X = 2\nu\{\lambda_X, \lambda_Y\}.$$

In order to avoid technical difficulties in defining the star representation (see below), we will assume our covariant star products to satisfy the following condition.

**Definition 2.3** *Let $(M, \omega, \mathfrak{g}, \lambda)$ be a strongly hamiltonian system. Let $\star_\nu$ be a $\mathfrak{g}$−covariant star product on $(M, \omega)$. We say that $\star_\nu$ has the property (B) if there exists an integer $N \in \mathbb{N}$ such that one has*

$$(\lambda_X \star_\nu u) \bmod(\nu^N) = (\lambda_X \star_\nu u) \bmod(\nu^{N+n});$$
$$(u \star_\nu \lambda_X) \bmod(\nu^N) = (u \star_\nu \lambda_X) \bmod(\nu^{N+n}),$$

*for all $X \in \mathfrak{g}$ and $u \in C^\infty(M)$ and $n \in \mathbb{N}$.*

In other words the series $\lambda_X \star_\nu u$ and $u \star_\nu \lambda_X$ stop at order $N$ independently of $u$ in $C^\infty(M)$. The following definition then makes sense [8].

**Definition 2.4** *Let $(M, \omega, \mathfrak{g}, \lambda)$ be a strongly hamiltonian system and $\star_\nu$ be a $\mathfrak{g}$-covariant star product on $(M, \omega)$ satisfying the property (B). One defines the representation $\mathrm{ad}_{\star_\nu}$ of $\mathfrak{g}$ on $E_\nu$ by*

$$\mathrm{ad}_{\star_\nu}(X)a := \frac{1}{2\nu}(\lambda_X \star_\nu a - a \star_\nu \lambda_X). \tag{5}$$

**Definition 2.5** *Let $G$ be a connected Lie group with Lie algebra $\mathfrak{g}$. Let $(M, \omega, \mathfrak{g}, \lambda)$ be a strongly hamiltonian system. Let $\star_\nu$ be a $\mathfrak{g}$-covariant star product on $(M, \omega)$ satisfying the property (B). The associated star representation of $G$ (if it exists) is the representation $\pi$ of $G$ on $E_\nu$ such that*

$$d\pi = \mathrm{ad}_{\star_\nu}.$$

## 3 Holonomy reducible symmetric orbits

Let $G$ be a connected simple Lie group. Let us denote by $\mathfrak{g}$ its Lie algebra. Let $\mathcal{O}$ be an adjoint orbit of $G$ in $\mathfrak{g}$. Choose a base point $o$ in $\mathcal{O}$ and denote by $\mathfrak{h}$ the Lie algebra of its stabilizer in $\mathfrak{g}$. Since $\mathfrak{g}$ is simple the Killing form $\beta$ establishes an equivariant linear isomorphism between $\mathfrak{g}$ and its dual $\mathfrak{g}^\star$. Therefore, every adjoint orbit can be identified with a coadjoint one. We will denote by $\omega^{\mathcal{O}}$ the KKS symplectic structure on $\mathcal{O}$ (cf. Example 1).

**Definition 3.1** *An adjoint orbit $\mathcal{O}$ in $\mathfrak{g}$ is symmetric if there exists an involutive automorphism $\sigma$ of $\mathfrak{g}$ such that $\mathfrak{h} = \{X \in \mathfrak{g} | \sigma(X) = X\}$. In this case, we denote by $\mathfrak{g} = \mathfrak{h} \oplus \mathfrak{q}$ the decomposition of $\mathfrak{g}$ induced by $\sigma$.*

One has $[\mathfrak{h}, \mathfrak{q}] = \mathfrak{q}$ and $[\mathfrak{q}, \mathfrak{q}] = \mathfrak{h}$ and the following is classical.

**Proposition 3.1** *Let $\mathcal{O}$ be a symmetric adjoint orbit of a simple Lie group $G$ and $(\mathfrak{g}, \sigma, \Omega)$ be its associated symplectic symmetric triple. The following assumptions are equivalent.*

(i) *the center $\mathfrak{z}(\mathfrak{h})$ of $\mathfrak{h}$ contains a non-compact element.*

(ii) *The subspace $\mathfrak{q}$ splits into a direct sum $\mathfrak{q} = \mathfrak{l} \oplus \mathfrak{l}'$ of isomorphic $\mathfrak{h}$-modules. One has $[\mathfrak{l}, \mathfrak{l}] = 0, [\mathfrak{l}', \mathfrak{l}'] = 0$ and both $\mathfrak{l}$ and $\mathfrak{l}'$ are $\beta$-isotropic and $\Omega$-Lagrangian subspaces of $\mathfrak{q}$.*

Such an orbit is called *holonomy reducible* if $\mathfrak{h}$ acts reducibly on $\mathfrak{q}$.

**Proposition 3.2** *([5]) Let $\mathcal{O}$ be a holonomy reducible symmetric orbit in $\mathfrak{g}$. We define the map $\phi : \mathfrak{q} = \mathfrak{l} \oplus \mathfrak{l}' \to \mathcal{O}$ by*

$$\phi(l, l') := Ad(exp(l).exp(l')).o$$

*Then $\phi$ is a Darboux chart on $(\mathcal{O}, \omega^{\mathcal{O}})$. Precisely, one has $\phi^{\star}\omega^{\mathcal{O}} = \Omega$.*

We transport the infinitesimal action of $\mathfrak{g}$ on $\phi(\mathfrak{q}) \subset \mathcal{O}$ to an infinitesimal action of $\mathfrak{g}$ on $\mathfrak{q}$ via $\phi$. One then gets a homomorphism of Lie algebra

$$\mathfrak{g} \to \mathcal{X}(\mathfrak{q}) : X \to X^{\star}.$$

Setting
$$\lambda_A(x) := \beta(\phi(x), A) \qquad (x \in \mathfrak{q}, A \in \mathfrak{g}),$$

one obtains the strongly hamiltonian system (cf.Definition 2.1) $(\mathfrak{q}, \Omega, \mathfrak{g}, \lambda)$. The main property of the Darboux chart $\phi$ is

**Proposition 3.3** *[5, 8] The Moyal star product on the symplectic vector space $(\mathfrak{q}, \Omega)$ is $\mathfrak{g}$-covariant. Moreover the Moyal star product in this case satisfies the property (B).*

**Definition 3.2** *Let $\mathcal{S}'_2(\mathfrak{q})$ be the space of distributions on $\mathfrak{q}$ which are tempered in the $\mathfrak{l}'$-variables w.r.t. the Lebesgue measure $dl'$ on $\mathfrak{l}'$. On $\mathcal{S}'_2(\mathfrak{q})$, we consider the partial Fourier transform $\mathcal{F} : \mathcal{S}'_2(\mathfrak{q}) \to \mathcal{S}'_2(\tilde{\mathfrak{q}})$ which reads formally as*

$$(\mathcal{F}u)(l, \eta) := \int_{\mathfrak{l}'} e^{-i\Omega(\eta, l')} u(l, l') dl' \qquad i := \sqrt{-1}.$$

Here $\tilde{\mathfrak{q}} := \mathfrak{l} \oplus \mathfrak{l}$ i.e. we identify the dual space $\mathfrak{l}'^{\star}$ with $\mathfrak{l}$ by use of $\Omega$.
We will also adopt the notation $\hat{u} := \mathcal{F}u$

We define the $\mathbb{R}$-isomorphism

$$\tilde{\mathfrak{q}} \to \mathfrak{l}^{\mathbb{C}} : (l, \eta) \to z = l + i\nu\eta$$

where the parameter $\nu$ is now considered as being *real*. In $\mathfrak{q} = \mathfrak{l} \oplus \mathfrak{l}'$, we choose an $\Omega$-symplectic basis $\{L_a, L'_a; 1 \leq a \leq n\}$ where $L_a \in \mathfrak{l}$ and $L'_a \in \mathfrak{l}'$. On $\mathfrak{l}$, we define the coordinate system $x = \Omega(x, L'_a)L_a =: x^a L_a$.
If $h$ is an element of $\mathfrak{h}$, we define its trace as

$$\mathrm{spur}(h) := \Omega([h, L_a], L'_a).$$

In this setting, one has the holomorphic constant vector field on $\mathfrak{l}^{\mathbb{C}}$:

$$\partial_{z^a} := \frac{1}{2\nu}(\nu(L_a)_l - i(L_a)_\eta) \quad (1 \leq a \leq n).$$

**Definition 3.3** *Considering $\mathfrak{l}^{\mathbb{C}} \subset \mathfrak{g}^{\mathbb{C}}$, we define, for all $A \in \mathfrak{g}^{\mathbb{C}}$, the polynomials on $\mathfrak{l}^{\mathbb{C}}$ :*

$$h_A^{\mathbb{C}}(z) := A_{\mathfrak{h}^{\mathbb{C}}} + [A_{\mathfrak{l}'^{\mathbb{C}}}, z] \in \mathfrak{h}^{\mathbb{C}}$$

$$l_A^{\mathbb{C}}(z) := A_{\mathfrak{l}^{\mathbb{C}}} + [A_{\mathfrak{h}^{\mathbb{C}}}, z] + \tfrac{1}{2}[z, [z, A_{\mathfrak{l}'^{\mathbb{C}}}]] \in \mathfrak{l}^{\mathbb{C}}$$

*where $z \in \mathfrak{l}^{\mathbb{C}}$ and $A = A_{\mathfrak{h}^{\mathbb{C}}} + A_{\mathfrak{l}^{\mathbb{C}}} + A_{\mathfrak{l}'^{\mathbb{C}}}$ according to the decomposition $\mathfrak{g}^{\mathbb{C}} = \mathfrak{h}^{\mathbb{C}} \oplus \mathfrak{l}^{\mathbb{C}} \oplus \mathfrak{l}'^{\mathbb{C}}$.*

**Definition 3.4** *For all $A \in \mathfrak{g}^{\mathbb{C}}$, we define the holomorphic vector field $\mathcal{Z}_A^{(\nu)} \in \Gamma(T^{1,0}(\mathfrak{l}^{\mathbb{C}}))$ by*

$$(\mathcal{Z}_A^{(\nu)})_z.f := (l_A^{\mathbb{C}}(z))^a (\partial_{z^a}.f)(z)$$

*where $z \in \mathfrak{l}^{\mathbb{C}}, f \in C^\infty(\mathfrak{l}^{\mathbb{C}}, \mathbb{C})$ and where, for all $w = w_1 + iw_2 \in \mathfrak{l}^{\mathbb{C}} = \mathfrak{l} \oplus i\mathfrak{l}$, we set $w^a := w_1^a + iw_2^a \quad (1 \le a \le n)$.*
*In the same way, extending the Killing form $\beta$ and the trace spur $\mathbb{C}$-linearly to $\mathfrak{h}^{\mathbb{C}}$, we define the complex polynomial of degree 1 :*

$$\tau_A^{(\nu)} := \frac{1}{2\nu}(\beta(h_A^{\mathbb{C}}, o) + \nu \, spur(h_A^{\mathbb{C}})).$$

**Proposition 3.4** *[5, 8] For all $A \in \mathfrak{g}$, one has*

$$\frac{1}{2\nu}\mathcal{F}(\lambda_A \star_\nu u) = \tau_A^{(\nu)}.\hat{u} + \mathcal{Z}_A^{(\nu)}.\hat{u}$$

*where $u$ is chosen in such a way that the LHS and the RHS make sense (e.g. $u \in \mathcal{S}(\mathfrak{q})$ the Schwartz space on $\mathfrak{q}$).*

We obtain a similar expression for the right regular $\star$ action.

## 4 Three different models for $\mathrm{ad}_{\star_\nu}$

From Proposition 3.4, straightforward computations yields

**Proposition 4.1** *For all $X \in \mathfrak{g}$ and $u \in C^\infty(\mathfrak{q})$ one has*

$$\begin{aligned}\mathrm{ad}_{\star_\nu}(X)u &= [(\Omega(L_a, X_{\mathfrak{l}'}) - \Omega([h_X, L_a], .))L'_a u + \Omega(l_X, L'_a)L_a u \\ &+ \nu^2 \Omega([[X_{\mathfrak{l}'}, L_a], L_b], L'_c)L'_a L'_b L_c u.\end{aligned}$$

*Or equivalently,*

$$\mathcal{F}(\mathrm{ad}_{\star_\nu}(X)u)(l, \eta) = \{i\beta([X_{\mathfrak{l}'}, \eta], o) + \mathrm{tr}(h_X(l))\}\hat{u}(l, \eta) \qquad (6)$$

$$+ \ \Omega(l_X(l) - \frac{\nu^2}{2}[[X_{\mathfrak{l}'}, \eta], \eta], L'_a)(L_a)_l.\hat{u}(l, \eta) + \Omega([h_X(l), \eta], L'_a)(L_a)_\eta.\hat{u}(l, \eta),$$

*where $(L_\eta.f)(l, \eta) := \frac{d}{dt}|_0 f(l, \eta + tL)$ and $1 \le a, b, c \le n$.*

Let us now turn to the particular case $G = SL(2,\mathbb{R})$. As usual we denote by $E, F$ and $H$ the generators of the Lie algebra $\mathfrak{sl}(2,\mathbb{R})$ satisfying the following commutation relations:

$$[E, H] = -2E, \quad [F, H] = 2F, \quad [E, F] = H.$$

We identify the vector spaces $\mathfrak{l}$ and $\mathfrak{l}'$ with $\mathbb{R} \cdot E$ and $\mathbb{R} \cdot F$ respectively and notice that the sub-algebra $\mathfrak{h}$ is isomorphic to $\mathbb{R} \cdot H$.

Consider the Lie algebra representation $\hat{\rho}_\nu$ of $\mathfrak{sl}(2,\mathbb{R})$ on the Schwartz space $\mathcal{S}(\mathfrak{l} \oplus \mathfrak{l})$ given by

$$\hat{\rho}_\nu(X)\hat{f}(l,\eta) := \mathcal{F}(\mathrm{ad}_{\star_\nu}(X)f). \tag{7}$$

According to (6) we have

$$\hat{\rho}_\nu(X)\hat{f}(l,\eta) = m_\nu(X)\hat{f}(l,\eta) + \mathcal{Y}_\nu(X)\hat{f}(l,\eta), \tag{8}$$

where

$$\begin{aligned}
m_\nu(X)f &:= (i\beta([X_{\mathfrak{l}'},\eta],o) + \mathrm{tr}(h_X(l)))f, \\
\mathcal{Y}_\nu(X)f &:= (\Omega(l_X(l) - \frac{\nu^2}{2}[[X_{\mathfrak{l}'},\eta],\eta],L_a')\partial_l + \Omega([h_X(l),\eta],L_a')\partial_\eta)f.
\end{aligned}$$

By identification of elements $l \cdot E$ and $\eta \cdot F$ with their coefficients $l$ and $\eta$ respectively we get

$$\mathcal{Y}_\nu(E) = \partial_l, \quad \mathcal{Y}_\nu(H) = 2l\partial_l + 2\eta\partial_\eta, \quad \mathcal{Y}_\nu(F) = (\nu^2\eta^2 - l^2)\partial_l - 2l\eta\partial_\eta.$$

Using the identification of $\mathfrak{l} \oplus \mathfrak{l}$ with $\mathfrak{l}^\mathbb{C}$ given by the $\mathbb{R}$-isomorphism $(l,\eta) \mapsto z = l + i\nu\eta$ we may write the ODE's describing the local flows of $\mathcal{Y}_\nu(X)$ as follows:

$$\text{For} \quad E \quad \dot{z} = 1, \quad \text{for} \quad H \quad \dot{z} = 2z, \quad \text{for} \quad F \quad \dot{z} = -z^2.$$

One recognizes the generating relations for three subgroups of conformal transformations of the upper half plane $\mathbb{H}^+$.

Under the same assumptions the multiplicative cocycle $m_\nu(X)$ is given by

$$m_\nu(X) = 2X_\mathfrak{h} - 2X_{\mathfrak{l}'}\mathrm{Re}(z) + \frac{i}{\nu}X_{\mathfrak{l}'}\mathrm{Im}(z).$$

Let us notice that the representation $\rho_\nu$ "looks like" a representation given by the parabolic induction. More precisely, consider the principal series representation of the complex Lie group $SL(2,\mathbb{C})$ acting on smooth functions by

$$\mathcal{P}^{\mu,\ell}(g)f(z) = |cz+d|^\mu \left(\frac{cz+d}{|cz+d|}\right)^\ell f\left(\frac{az+b}{cz+d}\right), \quad g^{-1} = \begin{pmatrix} a & b \\ c & d \end{pmatrix},$$

where $\mu \in \mathbb{C}, \ell \in \mathbb{Z}$. This action is unitary for $\mathrm{Re}(\mu) = -2$ (see [15] p.33). The corresponding infinitesimal action of the generator $F$ is given by

$$d\mathcal{P}^{\mu,\ell}(F)f(z) = (\mu\mathrm{Re}(z) + i\ell\mathrm{Im}(z))f(z) - z^2 f'(z).$$

So we have

**Theorem 4.2** *The star representation* $\mathrm{ad}_{\star_\nu}$ *exponentiates at the Lie group level and it is equivalent to the restriction of the principal series representation* $\mathcal{P}^{\mu,\ell}$ *of the Lie group* $G^{\mathbb{C}} = SL(2,\mathbb{C})$ *to its real non-compact form* $SL(2,\mathbb{R})$ *for* $\mu = -2, \ell = \nu^{-1}$. *The intertwining operator is the partial Fourier transform.*

This result can be generalized in the framework of so-called "conformal" or Kantor-Koecher-Tits Lie algebras. Indeed, such a Lie algebra is a symmetric graded Lie algebra

$$\mathfrak{g} = \mathfrak{g}_{-1} \oplus \mathfrak{g}_0 \oplus \mathfrak{g}_1,$$

with $[\mathfrak{g}_k, \mathfrak{g}_l] = \mathfrak{g}_{k+l}$ for $k \neq l \in \{-1, 0, 1\}$. Moreover, a Kantor-Koecher-Tits algebra has a canonical $\mathfrak{sl}_2$-triplet $(e, f, g), e \in \mathfrak{g}_{-1}, h \in \mathfrak{g}_0, f \in \mathfrak{g}_1$ defining this graduation:

$$\mathfrak{g}_k = \{X \in \mathfrak{g} | \, [h, X] = 2kX\}.$$

The isomorphic $\mathfrak{g}_0$-modules $\mathfrak{g}_{\pm 1}$ are endowed with the additional structure of a Jordan algebra. If such an algebra is euclidian (it carries a positive-definite bilinear form) then the real connected Lie group $G$ such that $\mathfrak{g} = \mathrm{Lie}(G)$ is the Lie group of holomorphic automorphisms of the associated hermitian symmetric space of tube type $G/K \simeq \mathfrak{g}_{-1} + iC$ where $C$ denotes the cone of positive elements in $\mathfrak{g}_{-1}$ (see [8, 11] for more details). Such a "conformal" Lie group is generated by translations of $\mathfrak{g}_{-1}$, linear transformations preserving the cone $C$ (the structure group) and the Jordan algebra inversion.

Using these three families of generators we can show that the star representation $\mathrm{ad}_{\star_\nu}$ of the Lie group of holomorphic transformations of an hermitian symmetric space of tube type is equivalent to the maximal principal series representation $\mathcal{P}_{\mu,\lambda}$ of the complex group $G^{\mathbb{C}}$ restricted to $G$ for some particular values of parameters $\mu$ and $\lambda$. Complete description of representations $\mathcal{P}_{\mu,\lambda}$ of conformal Lie groups in terms of corresponding Jordan algebras is done in [19].

Unitary irreducible representations of the maximal principal series of the complex Lie group $SL(2,\mathbb{C})$ become deeply reducible when restricted to the real subgroup $SL(2,\mathbb{R})$. We shall discuss the corresponding brunching laws. Moreover let us notice that an algebraic approach to the problem of restrictions of principal series representations of complex semi-simple Lie groups to their real forms is given in [17].

Back to the $SL(2,\mathbb{R})$-case, since $SU(1,1)$ and $SL(2,\mathbb{R})$ are conjugate in $SL(2,\mathbb{C})$ we can restrict all our considerations to $G = SU(1,1)$. Let $K = U(1)$ be the maximal compact subgroup of $G$ and let $\Gamma(G/K, \ell)$ denote the set of functions $f \in \Gamma(G)$ such that $f(gk) = \chi_\ell(k)f(g)$ where $\chi_\ell$ is a character of $K$ given by
$\chi_\ell \begin{pmatrix} e^{i\phi} & 0 \\ 0 & e^{-i\phi} \end{pmatrix} = e^{i\ell\phi}, \ell \in \mathbb{Z}, \phi \in [0, 2\pi]$.

**Theorem 4.3** *For* $\mathrm{Re}\,\mu = -2$ *the restriction of the principal series representation* $\mathcal{P}^{\mu,\ell}$ *of* $G^{\mathbb{C}}$ *to* $G$ *is equivalent to*

$$\mathcal{P}^{\mu,\ell}|_G \simeq L^2(G/K, \ell) \oplus L^2(G/K, -\ell).$$

For the proof, recall the construction of the parabolic induction for representations $\mathcal{P}^{\mu,\ell}$ (see for example [10]) in the compact picture. Let $P^{\mathbb{C}}$ be the Borel subgroup of $G^{\mathbb{C}}$ consisting of the lower triangular matrices $p = \begin{pmatrix} a & 0 \\ c & a^{-1} \end{pmatrix}$, with $c \in \mathbb{C}, a \in \mathbb{C}^*$.

The group $G^{\mathbb{C}}$ acts on $S = \{s \in \mathbb{C}^2 : \|s\| = 1\}$ and on $\tilde{S} = S/\sim$ (where $s \sim s'$ if there exists $\lambda \in \mathbb{C}, |\lambda| = 1$ such that $s = \lambda s'$) by

$$g \cdot s = g(s)/\|g(s)\|.$$

The action is transitive and, clearly, the stabilizer of the first basis vector $\tilde{e}_1$ is the maximal parabolic $P^{\mathbb{C}}$. So $\tilde{S} \simeq G^{\mathbb{C}}/P^{\mathbb{C}}$.
If $ds$ denotes the induced measure on $S$, then

$$d(g \cdot s) = \|g(s)\|^{-4} ds.$$

For $\mu \in \mathbb{C}$ and $\ell \in \mathbb{Z}$, let us define the character $\omega_{\mu,\ell}$ of $P^{\mathbb{C}}$ by

$$\omega_{\mu,\ell}(p) = |a|^{\mu}(a/|a|)^{\ell}.$$

Consider $\mathcal{P}^{\mu,\ell} = \mathrm{Ind}_{P^{\mathbb{C}}}^{G^{\mathbb{C}}} \omega_{\mu,\ell}$. Let $C^{\infty}(\tilde{S}, \ell)$ be the space of smooth functions $\varphi$ on $S$ satisfying

$$\varphi(\lambda s) = \lambda^{-\ell} \varphi(s), \tag{9}$$

for $s \in S, \lambda \in \mathbb{C}, |\lambda| = 1$. Then the representation $\mathcal{P}^{\mu,\ell}$ acts on it by

$$\mathcal{P}^{\mu,\ell}(g)\varphi(s) = \varphi(g^{-1} \cdot s)\|g^{-1}(s)\|^{\mu}.$$

Let $(\ ,\ )$ be the standard inner product on $L^2(S)$: $(\varphi_1, \varphi_2) = \int_S \varphi_1(s)\overline{\varphi_2(s)} ds$. This form is invariant with respect to $(\mathcal{P}^{\mu,\ell}, \mathcal{P}^{-\bar{\mu}-4,\ell})$. Therefore, if $\mathrm{Re}\,\mu = -2$, then the representation $\mathcal{P}^{\mu,\ell}$ is unitary.
The real form $G$ acts on $S$ too. Let $[\,,\,]$ be the sesqui-linear form on $\mathbb{C}^2$ given by $[z,w] = z_1\bar{w}_1 - z_2\bar{w}_2$. There are 3 orbits for this action given by

$$O_1 : [s,s] > 0 \qquad O_2 : [s,s] < 0 \qquad O_3 : [s,s] = 0.$$

We can identify the two first open orbits $O_1$ and $O_2$ with the group $G$ itself by $g \to g \cdot e_1$ and $g \to g \cdot e_2$ respectively. Notice also that for $k = \begin{pmatrix} e^{i\phi} & 0 \\ 0 & e^{-i\phi} \end{pmatrix} \in K$ one has $gk \cdot e_1 = e^{i\phi} g \cdot e_1$ and $gk \cdot e_2 = e^{-i\phi} g \cdot e_2$. The third orbit is closed and its measure is zero.
Consider a smooth, compactly supported function $\varphi$ on $O_1$ satisfying Relation (9). Then the function

$$A(\varphi) = \psi(s) := [s,s]^{-\mu/2} \varphi(s)$$

again satisfies Relation (9) and moreover

$$\psi(g^{-1} \cdot s) = \mathcal{P}^{\mu,\ell}(g)\varphi(s)[s,s]^{-\mu/2}, \qquad (g \in G).$$

Therefore the correspondence $A : \varphi \to \psi$ is a linear intertwining map from $\mathcal{P}^{\mu,\ell}|_G$ (on $O_1$) to the left regular action of the group $G$ on the space $C_c^{\infty}(G/K, \ell)$ which is dense in $L^2(G/K, \ell)$. Moreover, keeping in mind that $\frac{ds}{\|[s,s]\|^2}$ is a $G$-invariant measure on $\tilde{O}_1$, we have

$$\int_{\tilde{O}_1} |\psi(s)|^2 ds = \int_{\tilde{O}_1} [s,s]^{-\mathrm{Re}\,\mu-2} |\varphi(s)|^2 ds.$$

So, in the unitary case $\operatorname{Re}\mu = -2$, the map $A$ becomes a unitary intertwining operator. The same consideration for the second open orbit $O_2$ achieves the proof.

To prove the last theorem in a general setting of an arbitrary Hermitian symmetric space $G/K$ of tube type we need an explicit description of all open $G$-orbits on the complex flag manifold $G^{\mathbb{C}}/P^{\mathbb{C}}$. The answer to this question is given by the Theorem 5.9 in [26]. Assume that $Z = G^{\mathbb{C}}/P^{\mathbb{C}} \simeq G_{comp}/K$ is an irreducible hermitian symmetric spaces of compact type dual to $G/K$ ($G_{comp}$ denotes a compact real form of the connected, simply connected, complex, simple Lie group $G^{\mathbb{C}}$). Fix a Cartan involution $\theta$ corresponding to the decomposition $\mathfrak{g} = \mathfrak{k} \oplus \mathfrak{p}$, with $\mathfrak{k} = \operatorname{Lie}(K)$. There is a compact Cartan subalgebra $\mathfrak{t} \subset \mathfrak{k}$ of $\mathfrak{g}$. We denote $\Sigma(\mathfrak{g}^{\mathbb{C}}, \mathfrak{t}^{\mathbb{C}})$ the corresponding root system. According to Kostant's "cascade construction" one introduces the maximal set of strongly orthogonal non-compact positive roots in $\Sigma(\mathfrak{g}^{\mathbb{C}}, \mathfrak{t}^{\mathbb{C}})$. This set of cardinality $l = rank_{\mathbb{R}}\mathfrak{g}$ is given by $\Xi = \{\xi_1, \ldots, \xi_l\}$ where $\xi_1$ is the maximal root and $\xi_{m+1}$ is a maximal non-compact positive root orthogonal to $\{\xi, \ldots, \xi_m\}$. Any set of strongly orthogonal non-compact positive long roots in $\Sigma(\mathfrak{g}^{\mathbb{C}}, \mathfrak{t}^{\mathbb{C}})$ is $W(G,T)$-conjugate to a subset of $\Xi$, (where $W(G,T)$ denotes the Weyl group of the root system $\Sigma(\mathfrak{g}^{\mathbb{C}}, \mathfrak{t}^{\mathbb{C}})$). Let $z_0$ be the base point of the flag manifold $Z$, and let $c_\Gamma$ denotes the partial Cayley transform associated to a set of non-compact positive roots $\Gamma$.

**Theorem 4.4** *([26]) The $G$-orbits on $Z$ are just the orbits*

$$O_{\Gamma, \Delta} = G(c_\Gamma c_\Delta^2 z_0)$$

*where $\Gamma$ and $\Delta$ are disjoint subsets of $\Xi$. Two such orbits $O_{\Gamma,\Delta}$ and $O_{\Gamma',\Delta'}$ are equal if and only if cardinalities $|\Gamma| = |\Gamma'|$ and $|\Delta| = |\Delta'|$. An orbit $O_{\Gamma,\Delta}$ is open if and only if $\Gamma$ is empty. Any open orbit is maximal-dimensional.*

Therefore, according to theorems 4.3 and 4.4, there exists a finite number of irreducible representations $(\tau, V_\tau)$ of the maximal compact subgroup $K$ such that the restriction of unitary principal series representations of the complex hermitian group $G^{\mathbb{C}}$ to its non-compact real form decomposes into a direct sum of regular representations of $G$ on spaces of the type $L^2(G/K, \chi_\tau)$ of $L^2$-sections of the vector bundle $G \times_K V_\tau$.

Summarizing we can say that the star representation $\hat{\rho}_\nu$ with $(\nu = 1/\ell)$ yields an interpolation between the infinitesimal action of $\mathfrak{g}$ on the space of functions on the $H$-reducible coadjoint orbit (more precisely on the image of the Darboux chart) on one hand ($\nu = 0$) and the infinitesimal action of $\mathfrak{g}$ on space of functions on the hermitian symmetric space of tube type $G/K$ ($\ell = 0$) on the other hand.

Notice that a similar phenomenon occurs in the context of Berezin quantization. It was recently pointed out by Y.Neretin [18] for the hermitian symmetric space $X = U(p,q)/U(p) \times U(q)$ (which becomes of tube type for $p = q$). He considered the so-called canonical representations of the group $U(p,q)$ acting in Hilbert spaces of compactly supported smooth functions on $X$ endowed with scalar products given by positive-definite Berezin kernels $B_\lambda$. Notice that $B_\lambda$ is positive-definite precisely for $\lambda = \lambda_n$ negative integer and for $\lambda > Const > 0$ ($\lambda$ belonging to the Wallach set). Taking limits for $\lambda_n \to -\infty$ and $\lambda \to +\infty$ Neretin showed that these unitary representations give rise to an interpolation between regular representations of $U(p,q)$ and its dual compact real form $U(p+q)$ on $L^2(U(p,q)/U(p) \times U(q))$ and $L^2(U(p+q)/U(p) \times U(q))$ respectively.

## 5 The Poincaré disc

The decomposition of $L^2(G/K, \ell)$ into a sum of irreducible subspaces is explicitly described in [22] ( the more general case of the Plancherel formula for spherical functions on a Hermitian Lie group is studied in [21]).
Recall the following result by Takahashi (cf [22] 4.20-4.21). Let $f$ belongs to $L^2(G/K, \ell)$. Then

$$\|f_\ell\|^2 = \frac{1}{2\pi} \int_0^\infty \zeta_{\ell,1/2+i\tau}(f_\ell)\tau \tanh(\pi\tau)d\tau + \sum_{1 \leq p \leq |\ell|} (p-1/2)\zeta_{\ell,p}(f_\ell), \qquad (10)$$

where $\zeta_{n,s}$ denotes the positive definite spherical function on $G$ given by $\zeta_{n,s}(a_t) = (1-X)^s {}_2F_1(s+n, s-n, 1, X)$, where $X = \tanh^2(t/2)$ and $a_t = \begin{pmatrix} \cosh t & \sinh t \\ \sinh t & \cosh t \end{pmatrix}$.
It acts on $L^2(G/K, \ell)$ in sense of distributions.

Formula (10) shows that there exists an embedding $L^2(G/K) \hookrightarrow L^2(G/K, \ell)$ intertwining the corresponding regular representations. E.P. van den Ban told us an elegant way to realize any of such embedding in terms of Fourier transforms, we reproduce this here. For the details, we refer to [2] and literature cited there.

Let $G$ be a Hermitian simple Lie group. Write $G = KAN$ for an Iwasawa decomposition and consider $P = MAN$ the associated parabolic subgroup. Denote by $\mathfrak{a}$ the Lie algebra of $A$.
First, for $\lambda \in \mathfrak{a}_\mathbb{C}^*$, denote by $H_{1,\lambda}$ the space of measurable functions $f : G \to \mathbb{C}$ satisfying the following properties

(i) $f(xman) = a^{\lambda-\rho}f(x), \qquad (x \in G, (m,a,n) \in M \times A \times N)$,

(ii) $f|_K \in L^2(K)$.

This space is isomorphic to the Hilbert space $\mathcal{H} := L^2(K/M)$. The group $G$ unitarily acts on $H_{1,\lambda} = \mathcal{H}$ by left translations.
Second, consider an unitary irreducible representation $(V_\delta, \delta)$ of $K$. Consider the Hilbert space of sections $L^2(G/K, \delta)$ realized as equivariant functions

$$L^2(G/K, \delta) \subset \{f : G \to V_\delta | \; f(xk) = \delta(k)^{-1}f(x), \; x \in G, k \in K\}, \qquad (11)$$

and define
$$\mathcal{H}^\delta = \mathrm{Hom}_K(V_\delta, \mathcal{H}).$$

The latter is finite dimensional as isomorphic to $\mathrm{Hom}(V_\delta, 1)$ (by Frobenius reciprocity). The space $\mathcal{H}^\delta$ is non trivial if and only if $\delta$ has a non trivial $M$-fixed vector.
The operator-valued Fourier transform is then defined as a map $\mathcal{F}_\delta$ which associated to $f \in C_c^\infty(G/K, \delta)$ the meromorphic function $\mathcal{F}_\delta f : \mathfrak{a}_\mathbb{C}^* \to \mathrm{Hom}(\mathcal{H}^\delta, \mathcal{H})$ given by

$$\mathcal{F}_\delta f(\lambda)T = \int_G \pi_{1,\lambda}(x)T(f(x)) \; dx, \quad (T \in \mathcal{H}^\delta).$$

Observe that the transformation $\mathcal{F}_\delta$ intertwines the left action of $G$ on $C_c^\infty(G/K, \delta)$ with the representation $1 \otimes \pi_{1,\lambda}$ in $(\mathcal{H}^\delta)^* \otimes \mathcal{H}$. By the Plancherel formula [2], one sees that $\mathcal{F}_\delta$ actually extends to $L^2(G/K, \delta)$ as a partial isometry

$$L^2(G/K, \delta) \to L^2(i\mathfrak{a}^{*+}, \mathrm{Hom}(\mathcal{H}^\delta, \mathcal{H}), d\mu), \quad d\mu(\lambda) = \frac{d\lambda}{|c(\lambda)|^2},$$

(as usual $c(\lambda)$ denotes the Harish-Chandra **c**-function [15] p.279, and $\mathfrak{a}^{*+}$ the positive Weyl chamber) intertwining the left action and

$$\pi^\delta := \int_{i\mathfrak{a}^{*+}} (1^\delta \otimes \pi_{1,\lambda}) d\mu(\lambda).$$

We denote by $\mathcal{F}_\delta^*$ its transpose. The transformation $\mathcal{F}_1$ is even an isometry as it coincides with the usual Fourier transform.

Note that in the case $G = SL(2, \mathbb{R})$, the space $V_\delta$ is one-dimensional. Hence the Fourier transform $\mathcal{F}_\delta$ takes the simpler form

$$\mathcal{F}_\delta(f)(\lambda) = \int_G f(x) \pi_{1,\lambda}(x)|_{\mathcal{H}^\delta} \, dx.$$

**Theorem 5.1** *Let*
$$M^* : i\mathfrak{a}^{*+} \to \mathrm{Hom}(\mathcal{H}^\delta, \mathcal{H}^1)$$
*be a measurable function, almost everywhere bounded and such that $M^*(\lambda)$ is onto for almost every $\lambda$. Then, the linear map $L^2(G/K) \xrightarrow{T} L^2(G/K, \delta)$ given by*
$$Tf(\lambda) = \mathcal{F}_\delta^*(\mathcal{F}_1 f(\lambda) \circ M(\lambda)^*),$$
*is a continuous injective intertwiner. Moreover any such continuous injective intertwiner can be realized this way.*

We start by proving the second assertion. If $T : L^2(G/K) \to L^2(G/K, \delta)$ is a continuous intertwiner, then the map

$$\hat{T} : L^2(i\mathfrak{a}^{*+}, \mathrm{Hom}(\mathcal{H}^1, \mathcal{H}), d\mu) \to L^2(i\mathfrak{a}^{*+}, \mathrm{Hom}(\mathcal{H}^\delta, \mathcal{H}), d\mu)$$

defined as $\hat{T} := \mathcal{F}_\delta \circ T \circ \mathcal{F}^{-1}$ is a continuous linear map intertwining $\pi^1$ with $\pi^\delta$. Two representations $\pi_{1,\lambda}$ and $\pi_{1,\lambda'}$ ($\lambda \neq \lambda'$) being irreducible and nonequivalent, in vertue of the Schur lemma $\hat{T}$ is a multiplication operator. Hence, there exists a measurable function $M : i\mathfrak{a}^{*+} \to \mathrm{Hom}((\mathcal{H}^K)^*, (\mathcal{H}^\delta)^*)$ such that

$$(\hat{T}\varphi)(\lambda) = (M(\lambda) \otimes \mathrm{id}_\mathcal{H})(\varphi(\lambda))$$

viewing $\mathrm{Hom}(\mathcal{H}^\delta, \mathcal{H})$ as $(\mathcal{H}^\delta)^* \otimes \mathcal{H}$. Or equivalently,

$$\hat{T} = \varphi \circ M^*. \tag{12}$$

Note that since $\mathcal{H}^1$ is identified with the one dimensional space $\mathcal{H}^K$ of $K$-invariant elements in $\mathcal{H}$, one has that injectivity of $T$ implies injectivity of $M(\lambda)$ almost everywhere. We now prove the almost everywhere boundedness of $M$. Assume it is not the case, fix $\epsilon > 0$ and consider a subset $S$ of $i\mathfrak{a}^{*+}$ such that $0 < \mathrm{measure}(S) < \infty$ and $\|M^*(s)\| \geq \|\hat{T}\| + \epsilon$. Then one easily finds a unit element $\varphi \in L^2(i\mathfrak{a}^{*+}, \mathrm{Hom}(\mathcal{H}^1, \mathcal{H}), d\mu)$ such that $\|\varphi \circ M^*\|_2 \geq \|\hat{T}\| + \epsilon$ (choose $\varphi$ along the characteristic function of $S$). The latter contradicts equation (12). The rest of Theorem 5.1 being clear, we now may state the result announced in the introduction, namely

**Theorem 5.2** *Let $\ell \in \mathbb{Z}_0$ and denote by $(H(\mathbb{R}^2), \rho_\ell)$ the unitary representation of $\mathfrak{g} = \mathfrak{sl}(2,\mathbb{R})$ obtained by assigning the value $\frac{1}{\ell}$ to the parameter $\nu$ in the formal expression of the star representation $\mathrm{ad}_{\star_\nu}$ of $\mathfrak{g}$ in the Weyl-Moyal algebra (cf. Proposition 4.1, Theorems 4.2 and 4.3). Then, Theorems 4.2, 4.3 and 5.1 explicitly describe the injective continuous intertwiners*

$$L^2(G/K) \to (H(\mathbb{R}^2), \rho_\ell).$$